\documentclass[a4paper,12pt]{article}
\usepackage[utf8x]{inputenc}
\usepackage{amsmath}
\usepackage{amsfonts}
\usepackage{amssymb}
\usepackage{url}
\usepackage{graphicx}
\usepackage[export]{adjustbox}
\usepackage{subcaption}

\begin{document}

\title{Automatized Evaluation of Formalization Exercises in Mathematics}
\author{Merlin Carl}

\maketitle

\section{Introduction}

Learning the correct use of mathematical language frequently poses a challenge for beginner students. At the same time, it is a basic skill, required both for understanding mathematical texts and for presenting one's own work. 

In mathematical lectures and typical textbooks, this is rarely explictly discusses, though some offer a brief discussion, along with some formalization exercises (see, e.g., \cite{Ha}).



In this note, we present two pieces of software that pursue the goal to support beginner students in learning the use of formal language. 

The first one, called ``Math Dictations'' (a word that we learned from M. Junk, who used the concept (but no automatization thereof) in introductory courses at the university of Konstanz), challenges students to translate a proposition given in natural language, such as ``the real function $f$ is strictly increasing" into a quantifier formula such as $\forall{x}\forall{y}(x<y\rightarrow f(x)<f(y))$. It is similar to the formalization exercises that form part of the ``Mathematical Logic Tutor'' by A. Moreno (see \cite{BM}), but goes beyond this in (i) allowing first-order logic rather than propositional logic and (ii) using a restricted automated theorem prover for evaluating solutions, so that many solutions, rather than a single one, are recognized as correct answers. After the ``Math Dictations'' had been implemented and the first version of this article had been posted, we were made aware of the fact that this kind of formalization exercise is available in the Edukera system\footnote{See \url{https://www.edukera.com/}}. The Edukera formalization exercises work with an ATP in the background in full first-order logic and offer various contexts for exercises, among them also real functions with inequalities. Our description of ``Math Dictations'' should thus not be seen as a claim to priority of the concept, but   serves as an explanation of the program and in particular the syntax of the input language. However, we point out two differences between the ``Math Dictations'' and the formalization exercises in Edukera: First, the input in Edukera leaves little room for entering non well-formed formulas, while the ``Math Dictations'' allow a free input. Thus, in Edukera, there is more guidance for the user, while the ``Math Dictations'' offer more opportunities to make mistakes. Didactically, both approaches may well have complementary advantages and disadvantages. Second, while Edukera only returns a ``right-or-wrong''-feedback, the ``Math Dictations'' differentiate between (i) correct solutions, (ii) inputs that are necessary, but not sufficient, (iii) inputs that are sufficient, but not necessary and (iv) inputs that are neither necessary nor sufficient for the condition in question, which may help the user in improving a solution.

The second one, which, with a bow to the legacy of J. Conway and his ``Game of Life''\footnote{See, e.g., \url{https://bitstorm.org/gameoflife/}.} we call ``Game of Def'', has exercises that ask students to give descriptions of graphically depicted sets in a specified logical language with words such as ``right'', ``above'', ``neighbour'' or ``equal distance''.

Both programs are written in Prolog and form a part of the Diproche system, which is a proof checker for natural language proofs specifically adapted to the area of beginner exercises. The Diproche system is built by the example of the Naproche system due to P. Koepke, B. Schroeder, M. Cramer and others (see, e.g., \cite{Cr1} or \cite{CFKKSV}). The current Diproche version covers the topics of propositional calculus, Boolean set theory, sets and functions, elementary number theory, induction proofs and axiomatic geometry. Presentations of the checking mechanism and further components of Diproche can be found in \cite{CK} or \cite{C}.




\section{Math Dictations}



The idea of  ``math dictations'' is simple: The student is given a natural language expression, which she or he is then to translate it to a quantifier formula. The quantifier formula is then checked for correctness. As mentioned above, we first learned this concept from M. Junk in Konstanz. 

The automatization is rather straightforward: A dictation problem (Id,Nat,Formal,FreeVars) consists of an identifier Id, a natural language sentence (i.e., a string) Nat, a list of formal expressions in the internal Prolog list format Formal and a list of free variables that should occur in a solution. Here, Formal is a list of possible formalization of the sentence given in Nat. The reason we use a list rather than a single formalization is that we want to cover cases in which several substantially different approaches should count as equally correct. 

The accepted syntax of the current version is as follows:

\begin{itemize}
\item Small Latin letters are used for variables and constants; both variables and constants are terms.
\item Each natural number (written as a finite string of decimal digits) is a term. 
\item If $a$ and $b$ are terms and $a$ is not a number, then $a(b)$ is a term which describes the application of $a$ to $b$ (clearly, this only makes sense when $f$ is a function). 
\item When $a$ and $b$ are terms, then $a<b$, $a\leq b$, $a>b$, $a\geq b$ and $a=b$ are formulas. 
\item When $\phi$ and $\psi$ are formulas, then so are $(\phi\&\psi)$, $(\phi v\psi)$, $(\phi->\psi)$, $(\phi<->\psi)$ and $\sim\phi$.
\item When $\phi$ and $\psi$ are formulas and $x$ is a small Latin letter, then $Ax:\phi$ and $Ex:\phi$ are formulas.
\end{itemize}

All of these terms have their usualy meaning; as a convention, quantifiers range over real numbers. 
This language is sufficient to express, in the realm of real numbers, statement like the following:

\begin{enumerate}
\item Strictly between any two distinct real numbers, there is a third one.
\item $f$ is a strictly increasing function.
\item $f$ has a zero whenever $g$ has a zero. 
\item $f$ globally dominates $g$.
\item $f$ converges to $0$.
\end{enumerate}

Thus, this language is already sufficient for a variety of formalization exercises. 

In the program, the natural language formulation is displayed to the user, who also has a text window for entering a formula; clicking on the ``check'' button for the respective program, the checking is initiated and feedback is provided.

The checking works as follows: First, it is checked whether the input is a well-formed formula in which the right free variables appear (i.e., the same ones that appear in the natural language formulation). If not, an error message is displayed and no further processing takes place. Otherwise, the given expression $\phi$ is converted into an internal Prolog list format and a Prolog Tableau-prover\footnote{Unfortunately, the current version of the Tableau prover has a bug. It will be corrected soon.} (as, e.g., described in \cite{Fi}) is used to check, for each $\psi$ from $\psi_{1},...,\psi_{n}$ belonging to the list Formal in the specification of the problem, whether $\phi\rightarrow\psi$ and whether $\psi\rightarrow\phi$. If there are $i,j\in\{1,2,...,n\}$ such that both $\phi\rightarrow\psi_{i}$ and $\psi_{j}\rightarrow\phi$ can be verified, then the input is considered as correct and the user is congratulated for solving the problem. If there is $i\in\{1,2,...,n\}$ such that $\phi\rightarrow\psi_{i}$, but no $j\in\{1,2,...,n\}$ with $\psi_{j}\rightarrow\phi$, then a message is returned saying that $\phi$ is sufficient, but not necessary and that the input should be made more restrictive. 
If there is $i\in\{1,2,...,n\}$ such that $\psi_{i}\rightarrow\phi$, but no $j\in\{1,2,...,n\}$ with $\phi\rightarrow\psi_{j}$, then a message is returned saying that $\phi$ is necessary, but not sufficient and that the condition should be loosened. 
If there is neither such a $i$ nor such a $j$, the user is told that $\phi$ is neither sufficient nor necessary and that she or he should try again. 

Of course, the Tableau prover needs to be restricted in some way: First, due to the undecidablity of first-order logic, the checking might not terminate. Second, logical equivalence is a rather poor criterion for the adequacy of formalization. To take an extreme example, we should certainly not accept the statement of Fermat's last theorem as a formalization of example (1) claiming the density of the real numbers, just because both are provable! In our case, propositional equivalence is accepted without restriction, but the number of instantiations of universally quantified statements that can be used is restricted to $3$.\footnote{This value is not chosen for any particular reason, but experience so far shows that it is sufficient for all cases attempted so far and does not yield unacceptably long running times.}

\section{The `Game of Def'}

Math dictations as above only give a ```right'' or ``wrong'' answer, differentiated only by ``sufficient'' and ``necessary''. this is of little help in refining a wrong solution. 
it would be better if one could see what one actually defined, in contrast to what one was supposed to define. a good teacher could respond by giving examples that match the given solution but are not intended or that are wrongly not covered by an attempted formalization. however, automating this in general is quite difficult. For this reason, the ``Game of Def'' was designed. 

Different problem: Directly modelling a situation in a formal way that is not given by a natural language expression, but rather by a picture (or in some other way). 

The syntax of the formal language $\mathcal{L}_{\text{GD}}$ accepted by the system is as follows:

\begin{itemize}
\item Small latin letters denote variables and constants.
\item When $a$, $b$, $x$, $y$ are variables or constants, then rechts($a$,$b$), links($a$,$b$), ueber($a$,$b$), unter($a$,$b$), nachbar($a$,$b$) and dist($a$,$b$)=dist($x$,$y$) are formulas. (The meaning of these German terms will be explained below when we specify the semantics.)
\item When $\phi$ and $\psi$ are formulas, then $(\phi\&\psi)$, $(\phi v\psi)$, $(\phi->\psi)$, $(\phi<->\psi)$ and $\sim\phi$ are formulas.
\item When $\phi$ is a formula and $x$ is a small latin letter, then $Ex:\phi$ and $Ax:\phi$ are formulas. 
\end{itemize}

This syntax is adhered to strictly. No omission of brackets, e.g. by priority rules, or addition of extra brackets etc. are allowed. Though it would not be difficult to somewhat loosen those rule, this is in line with the didactical goal of helping to get used to expressing oneself within the borders of a formalism.\footnote{As it turns out, some of the advanced levels also raised the interest of advanced mathematicians, who took it as a kind of puzzle game. If this interest persists, loosening the syntactic rules will be reconsidered.} The somewhat odd notation for the existential and universal quantifier and the logical junctors is due to the implementation in Prolog. An improved interface with a more appealing input format is certainly desirable, though it should be kept in mind that beginners should not be expected to be familiar with LaTeX.

The semantics now works as follows: The domain on which the game is played is a $21\times21$-square grid $G$, with the middle marked with ``$u$''. Variables and constants refer to squares in this grid. Then:

\begin{itemize}
\item $a=b$ means that $a$ and $b$ denote the same square. 
\item rechts($a$,$b$) means that the square $b$ is somewhere to the right, but in the same row as, $a$; i.e., if one would use coordinates (which the game syntax does not), we would say that the $x$-coordinate of $b$ is larger than that of $a$, while the $y$-coordinates agree. 
\item links($a$,$b$) means that the square $b$ is somewhere to the left, but in the same row as, $a$.
\item ueber($a$,$b$) means that the square $b$ is somewhere above, but in the same column as, $a$.
\item unter($a$,$b$) means that the square $b$ is somewhere below, but in the same column as, $a$.
\item nachbar($a$,$b$) means that $a$ and $b$ are neighbours, i.e. share exactly one common border line. In coordinates, that means that they have one common coordinate, while they differ by $1$ in the other. 
\item dist($a$,$b$)=dist($x$,$y$) means that $a$ and $b$ lie in the same row or column, that $x$ and $y$ lie in the same row or column, and that the distance from $a$ to $b$ is the same as the distance from $x$ to $y$.
\end{itemize}

Junctors and quantifiers have their usual meaning; note that universal and exisential quantifiers only quantify over squares in the grid, not some infinite extension thereof. Thus, there are squares with no right neighbours etc. Formulas that contain more than $2$ nested quantifiers are accepted syntactically, but their semantic evaluation - which is based on an exhaustive search whenever nested quantifiers are involved - takes too long for all practical purposes. Thus, nesting more than two quantifiers should be avoided and is also not required for any solution.

The ``Game of Def'' now works as follows: In each exercise, one is given an image of the grid, with some squares marked yellow. Some of the squares may be labeled by letters, which means that those letters are constant letters that can be used as parameters. In addition, one is given an informal description of the set $Y$ of yellow squares in natural language (currently German). The task is then to write down a $\mathcal{L}_{\text{GD}}$-formula $\phi(x)$ with exactly one free variable $x$ (the choice of the variable is up to the user with the only restriction that constant letters used in the exercise description cannot be used) such that $\{x\in G:\phi(x)\}=Y$. 

Users can write a string into an input window and press the ``check'' button. If the input is not a $\mathcal{L}_{\text{GD}}$-formula or it does not have exactly one free variable, an error message is displayed and no further processing takes place. Otherwise, let us denote by $\phi(x)$ the input formula and by $U_{\phi}$ the set described by it. The system then does the following:

\begin{itemize}
\item Squares in $U_{\phi}\cap Y$ are colored green.
\item Squares in $U_{\phi}\setminus Y$ are colored red.
\item Squares in $Y\setminus U_{\phi}$ remain yellow.
\end{itemize}

Furthermore, the user receives the following text feedback:

\begin{itemize}
\item When $Y=U_{\phi}$, (s)he is congratulated that the solution is correct.
\item When $Y\subsetneq U_{\phi}$, a message is returned saying that the given condition is necessary, but not sufficient and that further restriction should be imposed.
\item When $U_{\phi}\subsetneq Y$, a message is returned saying that the given condition is sufficient, but not necessary and that it should be made more inclusive.
\item When none of the above cases hold, the user is told to try again. 
\end{itemize}

Here is an example of an exercise with the feedback as it is returned to the user: 

\includegraphics[angle=90, width=15cm, height=25cm ]{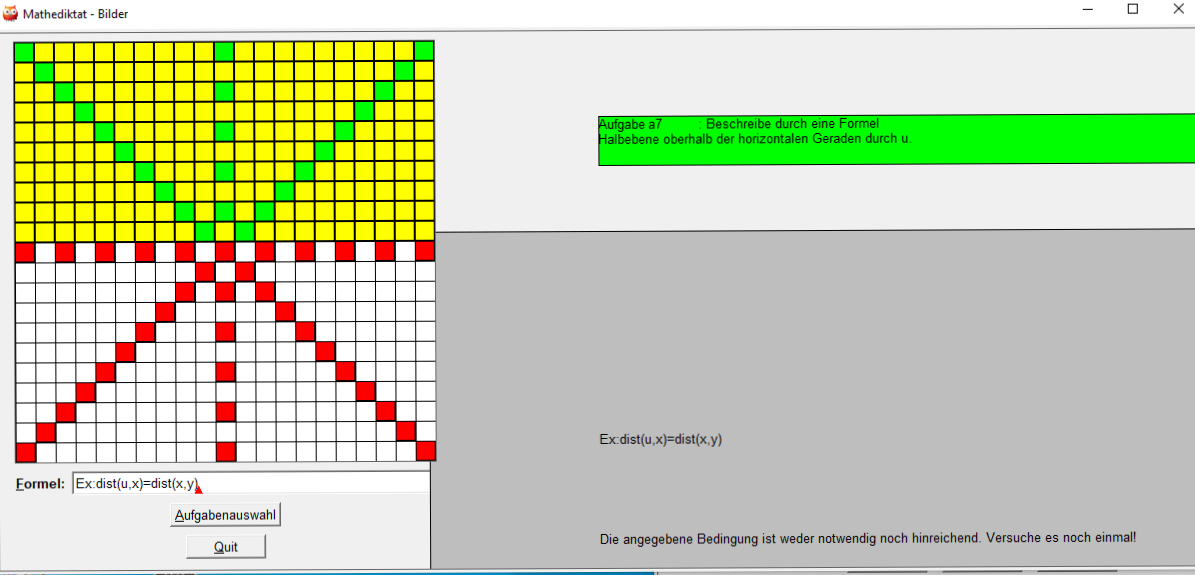}

The interested reader may now want to entertain her- or himself with the following exercises, which are part of the current version of the system:

\begin{figure}[h]

\begin{subfigure}{0.5\textwidth}
\includegraphics[width=0.9\linewidth, height=4cm]{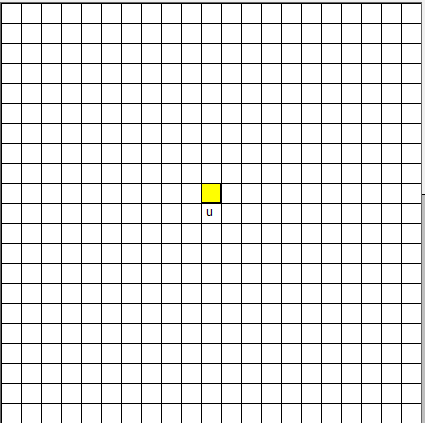} 
\caption{Problem1}
\label{fig:Problem 1}
\end{subfigure}
\begin{subfigure}{0.5\textwidth}
\includegraphics[width=0.9\linewidth, height=4cm]{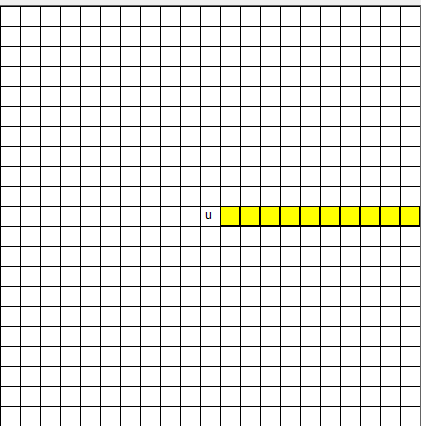}
\caption{Problem 2}
\label{fig:Problem 2}
\end{subfigure}

\end{figure}

\begin{figure}[h]

\begin{subfigure}{0.5\textwidth}
\includegraphics[width=0.9\linewidth, height=4cm]{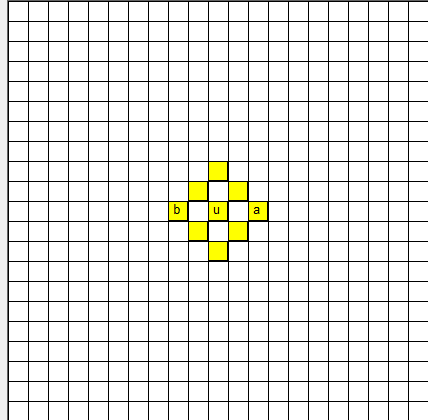} 
\caption{Problem3}
\label{fig:Problem 3}
\end{subfigure}
\begin{subfigure}{0.5\textwidth}
\includegraphics[width=0.9\linewidth, height=4cm]{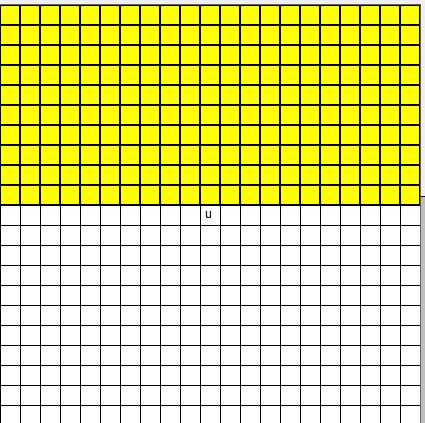}
\caption{Problem 4}
\label{fig:Problem 4}
\end{subfigure}

\end{figure}

\begin{figure}[h]

\begin{subfigure}{0.5\textwidth}
\includegraphics[width=0.9\linewidth, height=4cm]{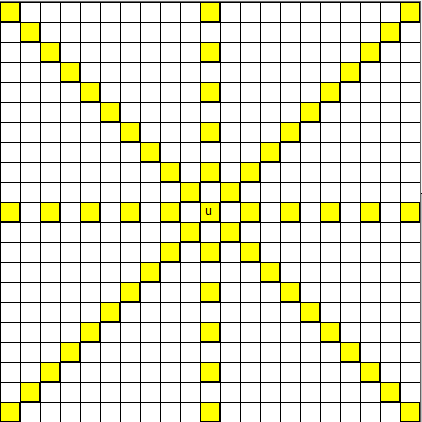} 
\caption{Problem5}
\label{fig:Problem 5}
\end{subfigure}
\begin{subfigure}{0.5\textwidth}
\includegraphics[width=0.9\linewidth, height=4cm]{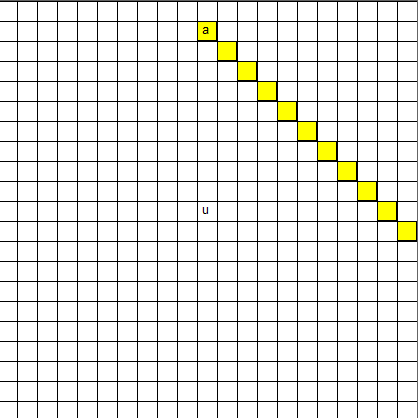}
\caption{Problem 6}
\label{fig:Problem 6}
\end{subfigure}

\end{figure}

\begin{figure}[h]

\begin{subfigure}{0.5\textwidth}
\includegraphics[width=0.9\linewidth, height=4cm]{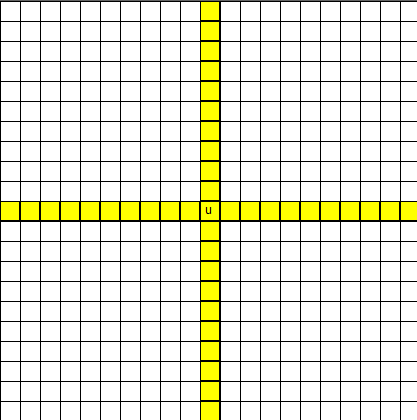} 
\caption{Problem7}
\label{fig:Problem 7}
\end{subfigure}
\begin{subfigure}{0.5\textwidth}
\includegraphics[width=0.9\linewidth, height=4cm]{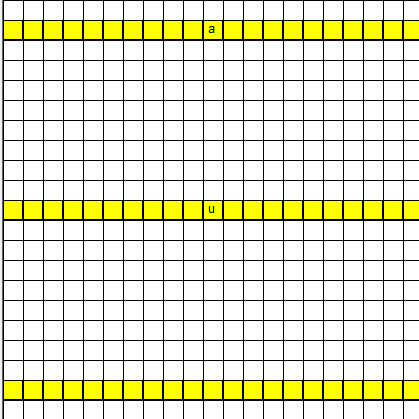}
\caption{Problem 8}
\label{fig:Problem 8}
\end{subfigure}

\end{figure}

\begin{figure}[h]

\begin{subfigure}{0.5\textwidth}
\includegraphics[width=0.9\linewidth, height=4cm]{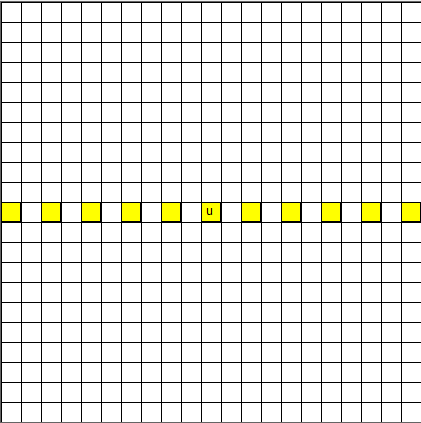} 
\caption{Problem9}
\label{fig:Problem 9}
\end{subfigure}
\begin{subfigure}{0.5\textwidth}
\includegraphics[width=0.9\linewidth, height=4cm]{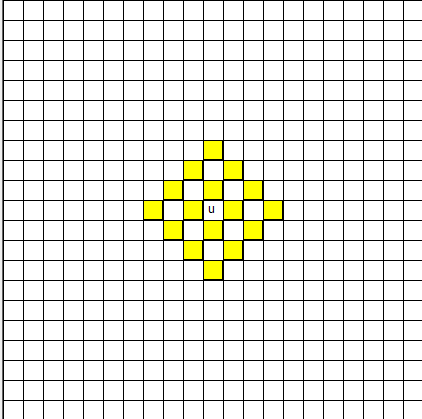}
\caption{Problem 10}
\label{fig:Problem 10}
\end{subfigure}

\end{figure}

\begin{figure}[h]

\begin{subfigure}{0.5\textwidth}
\includegraphics[width=0.9\linewidth, height=4cm]{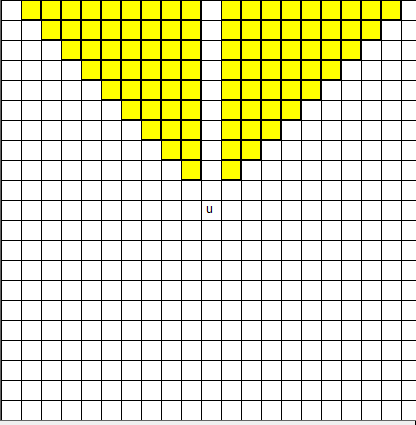} 
\caption{Problem11}
\label{fig:Problem 11}
\end{subfigure}
\begin{subfigure}{0.5\textwidth}
\includegraphics[width=0.9\linewidth, height=4cm]{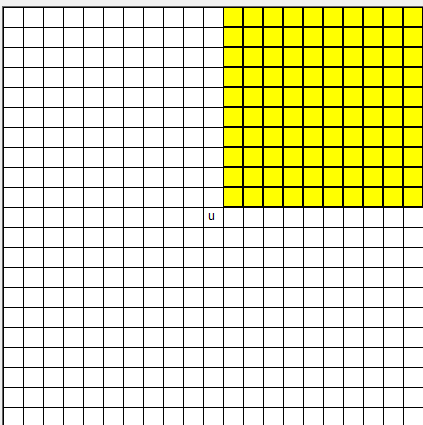}
\caption{Problem 12}
\label{fig:Problem 12}
\end{subfigure}

\end{figure}

\clearpage

\newpage

\section{Further Work}

Clearly, the possibilities of using automated theorem provers and truth predicate evaluation in supporting formalization exercises are endless. In particular, it is easy to extend the syntax of the math dictation program to comprise other areas of mathematics, like number theory or geometry. Concerning the Game of Def, it would be desirable to get rid of the limited number of nested quantifiers by improving the running time of the evaluation algorithm. 



There is a more general topic in the background here, which we plan to take up in future work: Namely, systematically look for theories that are both simple in terms of model theory and complexity theory ($o$-minimality, quantifier elimination and decidability (see, e.g., \cite{Ma}) seem to be particularly relevant properties) and didactically suitable in that their realm of objects is either known to or easy to explain to beginner students and that they allow for many non-trivial, but realistically solvable formalization exercises, preferable those with a visualizable aspect. The theories of Presburger arithmetic and real closed fields may be suitable candidates, provided that the complexity issues (Presburger arithmetic has a double-exponential lower time bound on a decision algorithm, see \cite{FR}; however, the situation is considerably less bad in the case of real closed fields, see, e.g., \cite{Gr}) turn out to be irrelevant for the intended application (simple formalization exercises). We hope for a stimulating interaction of mathematical logic (in particular model theory), computer science and the didactics of mathematics. 


\end{document}